\documentclass[12pt]{article}
\usepackage{amsmath}
\usepackage{amssymb}
\usepackage{latexsym}
\usepackage{amsthm}
\usepackage{mathrsfs}
\usepackage{enumitem}
\usepackage{mathdots}
\usepackage[colorlinks,
            linkcolor=blue,
            anchorcolor=blue,
            citecolor=blue
            ]{hyperref}

\newtheorem{thm}{Theorem}
\newtheorem{lem}[thm]{Lemma}

\newtheorem{cor}[thm]{Corollary}

\newcommand{\sep}{\preceq}

\begin{document}

\begin{center}
{\large \bf  On the real-rootedness of  the local $h$-polynomials of\\
edgewise subdivisions}
\end{center}

\begin{center}
Philip B. Zhang\\[6pt]

College of Mathematical Science \\
Tianjin Normal University, Tianjin  300387, P. R. China\\[8pt]

Email: {\tt zhang@tjnu.edu.cn}
\end{center}

\noindent\textbf{Abstract.}
Athanasiadis conjectured that, for every positive integer $r$, the local $h$-polynomial of the $r$th edgewise subdivision of any simplex has only real zeros.
In this paper, based on the theory of interlacing polynomials, we prove that a family of polynomials related to the desired local $h$-polynomial is interlacing and hence confirm Athanasiadis' conjecture.

\noindent \emph{AMS Classification 2010:} 26C10, 05E45, 05A15


\section{Introduction}

The objective of this paper is to prove a real-rootedness conjecture of Athanasiadis \cite{Athanasiadis2017survey} about the local $h$-polynomials of edgewise subdivisions of simplices.

Let us first review some background.
 Let {$\Delta$} be a $(d-1)$-dimensional simplicial complex with $f_i$ faces of dimension $i$ with  $f_{-1}=1$ by convention.
The \emph{$h$-polynomial} of $\Delta$ is defined as  $h(\Delta,x) = \sum_{i=0}^{d} f_{i-1} x^i (1-x)^{d-i}.$
The notion of local $h$-polynomials was introduced by Stanley \cite{Stanley1992Subdivisions}  in the study of the face enumeration of subdivisions of simplicial complexes.
Given an $n$-element set $V$, let $2^V$  be the abstract simplex consisting of all subsets of the set $V$ and
let $\Gamma$ be a simplicial subdivision of the simplex $2^V$.
The \emph{local $h$-polynomial} of $\Gamma$ is defined as
\begin{align*}
 \ell_V \left(\Gamma, x \right) = \sum_{F \subset V} (-1)^{n-|F|} \, h\left(\Gamma_{F},x\right),
\end{align*}
where $\Gamma_{F}$ is the restriction of $\Gamma$ to the face $F\in 2^{V}$.
In a recent survey, Athanasiadis \cite[Section 4]{Athanasiadis2017survey} studied several interesting examples of local $h$-polynomials and asked whether these polynomials have only real zeros.

 This paper is concerned with the local $h$-polynomials of the $r$th edgewise subdivision $(2^V)^{\langle r \rangle}$ of the simplex $2^V$.
Edgewise subdivision has appeared in several mathematical contexts, see \cite{Beck2010log, Brenti2009Veronese, Brun2005Subdivisions, Eisenbud1994Initial, Jochemko2018Real, Kubitzke2012Enumerative}.
One of its properties is that its faces $F$  are divided into $r^{\dim(F)}$  of the same dimension.
Athanasiadis \cite{Athanasiadis2014Edgewise, Athanasiadis2016local} showed that
\begin{align}\label{eq:localedge}
  \ell_V \left((2^V)^{\langle r \rangle}, x \right)=
  \ {\rm E}_r \left(\, (x + x^2 + \cdots + x^{r-1})^n \right),
\end{align}
where ${\rm E}_r$ is a linear operator defined on polynomials by setting ${\rm E}_r (x^n) = x^{n/r}$, if $r$ divides $n$,
and ${\rm E}_r (x^n) = 0$ otherwise.
The local $h$-polynomial $\ell_V \left((2^V)^{\langle r \rangle}, x \right)$ can also be interpreted combinatorially as the ascent generating function of certain Smirnov words \cite{Linusson2012Rees}, see \cite[Theorem 4.6]{Athanasiadis2017survey} for more details. We would like to point out that the real-rootedness of $\ {\rm E}_r \left(\, (1+x + x^2 + \cdots + x^{r-1})^n \right)$, which is slightly different from the right hand side of \eqref{eq:localedge}, has been studied in \cite{Jochemko2018Real,  ZhangInterlacing} and references therein.

The main result of this paper is as follows.
\begin{thm}\label{main-thm}
For any positive integer $r$, the local $h$-polynomial $\ell_V \left((2^V)^{\langle r \rangle}, x \right)$ has only real zeros.
\end{thm}
The above theorem provides an affirmative answer to Athanasiadis' conjecture. Note that this conjecture was also proved independently by Leander \cite{Leander2016Compatible}.

\section{Proof of  Theorem \ref{main-thm}}\label{sect-proof}
In this section, we shall prove Theorem \ref{main-thm}.
We first review the theory of interlacing polynomials, especially some useful criteria. We proceed to introduce a sequence of polynomials and prove that these polynomials are interlacing.
Finally, we  show that  one of these polynomials is just the polynomial ${\rm E}_r \, (x + x^2 + \cdots + x^{r-1})^n$.
Our proof of Theorem \ref{main-thm} is based on the theory of interlacing polynomials, which has been widely used to prove the real-rootedness of several polynomials arising in  combinatorics (\cite{Jochemko2018Real, Savage2015s, Yang2014Mutual, Zhang2015real}).

Given two real-rooted polynomials $f(x)$ and $g(x)$ with positive leading coefficients, let $\{u_i\}$ and $\{v_j\}$ be the set of zeros of $f(x)$ and $g(x)$, respectively.
We say that \emph{$g(x)$ interlaces $f(x)$}, denoted $g(x)\sep f(x)$, if either
$\deg f(x)=\deg g(x)=d$ and
\begin{align*}
v_d\le u_d \le v_{d-1}\le\cdots\le v_2\le u_2\le v_1\le u_1,
\end{align*}
or $\deg f(x)=\deg g(x)+1=d$ and
\begin{align*}
u_{d}\le v_{d-1}\le\cdots\le v_{2}\le u_{2}\le v_{1}\le u_{1}.
\end{align*}
For convention, we let $0\ll  f$ and $f\ll  0$ for any real-rooted polynomial $f$.
A sequence of real polynomials $(f_{1}(x),\dots,f_{m}(x))$
with positive leading coefficients is said to be  \emph{interlacing} if $f_{i}(x)\ll  f_{j}(x)$ for all $1\le i<j\le m$.

Following Br{\"a}nd{\'e}n \cite{Braenden2015Unimodality}, let $\mathcal{F}_n$ be the set of all interlacing
sequences $(f_i)_{i=1}^n$ of polynomials, and  $\mathcal{F}_n^+$ be the subset  of  $(f_i)_{i=1}^n \in \mathcal{F}_n$ such that  the coefficients of $f_i$ are all  nonnegative  for all $1 \leq i \leq n$.
A central problem in this area is to characterize  $m \times n$ matrices $G= (G_{ij}(x))$ of polynomials which maps  $\mathcal{F}_n^+$ to $\mathcal{F}_m^+$ via a matrix multiplication as follows:
$$
G\cdot (f_1, \ldots, f_n)^T = (g_1, \ldots, g_m)^T.
$$
Br{\"a}nd{\'e}n \cite[Theorem 8.5]{Braenden2015Unimodality} gave a characterization of the case when the polynomials considered have all nonnegative coefficients.

\begin{lem} [{\cite[Theorem 8.5]{Braenden2015Unimodality}}]\label{pb}
Let  $G= (G_{ij}(x))$  be an $m \times n$ matrix of polynomials. Then $G : \mathcal{F}_n^+ \rightarrow \mathcal{F}_m^+$  if and only if
\begin{enumerate}
\item $G_{ij}(x)$ has nonnegative coefficients for all $1\le i \le m$ and $1\le j \le n$, and
\item for all $\lambda, \mu >0$, $1\leq i< j \leq n$ and $1\leq k< \ell \leq m$,
\begin{equation}\label{pbxvw}
(\lambda x + \mu) G_{kj}(x) + G_{\ell j}(x)\ll  (\lambda x + \mu) G_{ki}(x) + G_{\ell i}(x).
\end{equation}
\end{enumerate}
\end{lem}

Given a polynomial $f(x)$, there exist uniquely determined polynomials  $f^{\langle r,0 \rangle}(x)$, $f^{\langle r,1 \rangle}(x)$, $ \ldots,$  $f^{\langle r,r-1 \rangle}(x)$ such that
\begin{align*}
	f(x) = f^{\langle r,0 \rangle}(x^r)+xf^{\langle r,1 \rangle}(x^r)+ \cdots +x^{r-1}f^{\langle r,r-1 \rangle}(x^r).
\end{align*}
In order to prove Theorem \ref{main-thm}, we next give the following result which plays a key role in our proof.
\begin{thm}\label{key}
Let $r$ and $\ell$ be two  positive integers with $\ell \le r-1$.
Suppose that $f(x)$ and
$g(x)$ are two polynomial  with nonnegative coefficients  satisfying
\begin{align}
 \left(1+x+\dots+x^{\ell}\right) f(x)  = g(x). \label{eq:fg}
\end{align}
If the sequence $\left(f^{\langle r, r-1 \rangle}(x), \dots, f^{\langle  r,1 \rangle}(x), f^{\langle r,0 \rangle}(x)\right)$ is  interlacing, then so is $\big( g^{\langle r, r-1 \rangle}(x)$, $\dots, g^{\langle  r,1 \rangle}(x), g^{\langle r,0 \rangle}(x) \big)$.
\end{thm}
\begin{proof}
For any $1\le i \le r$, taking all the terms of form  $x^{mr-i}$ where $m$  is a non-negative integer from both sides of  \eqref{eq:fg}, we get that 
\begin{align*}
	x^{r-i} g^{\langle r,r-i \rangle} (x^r)  & = \sum_{j=1}^{i+\ell-r} x^{r+j-i}  \cdot x^{r-j}  f^{\langle  r,r-j \rangle}(x^r)  +  \sum_{j=i}^{i+\ell}  x^{j-i} \cdot x^{r-j}   f^{\langle  r,r-j \rangle}(x^r)  \\[5pt]
	   & = x^{2r-i} \sum_{j=1}^{i+\ell-r} f^{\langle  r,r-j \rangle}(x^r)
	  + x^{r-i} \sum_{j=i}^{i+\ell}    f^{\langle  r,r-j \rangle}(x^r), 
\end{align*}
where $r+j-i\le \ell$ for $1 \le j\le i+\ell  -r $ in the first summation and  $j-i \le \ell$ for $i \le j\le i+\ell$ in the second summation, 
and hence 
\begin{align*}
 g^{\langle r,r-i \rangle} (x)   = x \sum_{j=1}^{i+\ell-r} f^{\langle  r,r-j \rangle}(x)
+ \sum_{j=i}^{i+\ell}    f^{\langle  r,r-j \rangle}(x).
\end{align*}
Thus  we obtain an alternative expression of \eqref{eq:fg} in a matrix form
\begin{align*}
G\cdot \left(f^{\langle r, r-1 \rangle}(x), \dots, f^{\langle  r,1 \rangle}(x), f^{\langle r,0 \rangle}(x)\right)^T
= \left( g^{\langle r, r-1 \rangle}(x), \dots, g^{\langle  r,1 \rangle}(x), g^{\langle r,0 \rangle}(x)\right)^T,
\end{align*}
where $G= \left( G_{i,j}(x) \right)$ is a square matrix of order $r$ with 
$$G_{i,j}(x)=\left\{  
	\begin{array}{cc}
		1, & i \le j \le i+\ell,  \\
		x, & j\le i+\ell-r, \\
		0, & \mbox{otherwise}. \\
	\end{array}
\right.$$
One example of $G$ for  $r=9$ and $\ell=5$ is as follows:
$$
\left(
\begin{array}{ccccccccc}
1 & 1 & 1 & 1 & 1 & 1 & 0 & 0 & 0 \\
0 & 1 & 1 & 1 & 1 & 1 & 1 & 0 & 0 \\
0 & 0 & 1 & 1 & 1 & 1 & 1 & 1 & 0 \\
0 & 0 & 0 & 1 & 1 & 1 & 1 & 1 & 1 \\
x & 0 & 0 & 0 & 1 & 1 & 1 & 1 & 1 \\
x & x & 0 & 0 & 0 & 1 & 1 & 1 & 1 \\
x & x & x & 0 & 0 & 0 & 1 & 1 & 1 \\
x & x & x & x & 0 & 0 & 0 & 1 & 1 \\
x & x & x & x & x & 0 & 0 & 0 & 1 \\
\end{array}
\right) .$$
%
All the possible two-by-two submatrices of $G$  are
$$
\left(
\begin{array}{cc}
	1 & 1 \\
	0 & 1 \\
\end{array}
\right),
\left(
\begin{array}{cc}
	1 & 0 \\
	0 & 1 \\
\end{array}
\right),
\left(
\begin{array}{cc}
	1 & 0 \\
	0 & 0 \\
\end{array}
\right),
\left(
\begin{array}{cc}
	1 & 1 \\
	1 & 1 \\
\end{array}
\right),
\left(
\begin{array}{cc}
	1 & 0 \\
	1 & 1 \\
\end{array}
\right),$$
$$
\left(
\begin{array}{cc}
	1 & 0 \\
	1 & 0 \\
\end{array}
\right),
\left(
\begin{array}{cc}
	0 & 0 \\
	1 & 0 \\
\end{array}
\right),
\left(
\begin{array}{cc}
	0 & 0 \\
	0 & 0 \\
\end{array}
\right),
\left(
\begin{array}{cc}
	1 & 1 \\
	0 & 0 \\
\end{array}
\right),
\left(
\begin{array}{cc}
	0 & 0 \\
	1 & 1 \\
\end{array}
\right),$$
$$
\left(
\begin{array}{cc}
0 & 0 \\
0 & 1 \\
\end{array}
\right),
\left(
\begin{array}{cc}
0 & 1 \\
0 & 0 \\
\end{array}
\right),
\left(
\begin{array}{cc}
0 & 1 \\
0 & 1 \\
\end{array}
\right),
\left(
\begin{array}{cc}
1 & 1 \\
x & 0 \\
\end{array}
\right),
\left(
\begin{array}{cc}
1 & 1 \\
x & 1 \\
\end{array}
\right),$$
$$\left(
\begin{array}{cc}
	1 & 0 \\
	x & 1 \\
\end{array}
\right),
\left(
\begin{array}{cc}
	1 & 0 \\
	x & 0 \\
\end{array}
\right),
\left(
\begin{array}{cc}
0 & 0 \\
x & 1 \\
\end{array}
\right), 
\left(
\begin{array}{cc}
0 & 0 \\
x & 0 \\
\end{array}
\right),
\left(
\begin{array}{cc}
0 & 1 \\
x & 0 \\
\end{array}
\right),$$
$$\left(
\begin{array}{cc}
0 & 1 \\
x & 1 \\
\end{array}
\right), 
\left(
\begin{array}{cc}
x & 0 \\
x & 0 \\
\end{array}
\right),
\left(
\begin{array}{cc}
x & 1 \\
x & 0 \\
\end{array}
\right),
\left(
\begin{array}{cc}
x & 1 \\
x & 1 \\
\end{array}
\right),
\left(
\begin{array}{cc}
1 & 1 \\
x & x \\
\end{array}
\right),$$
$$\left(
\begin{array}{cc}
	0 & 0 \\
	x & x \\
\end{array}
\right),
\left(
\begin{array}{cc}
0 & 1 \\
x & x \\
\end{array}
\right),
\left(
\begin{array}{cc}
	x & 0 \\
	x & x \\
\end{array}
\right),
\left(
\begin{array}{cc}
	x & 1 \\
	x & x \\
\end{array}
\right),\left(
\begin{array}{cc}
	x & x \\
	x & x \\
\end{array}
\right).$$
One can check all these submatrices satisfy the condition \eqref{pbxvw} of
Lemma \ref{pb}.
We take the first matrix as an example, and all the other matrices can be treated similarly.
We need to check for all $\lambda,\mu>0$ the following interlacing relation is satisfied,
	\begin{align*}
		\lambda x + \mu + 1\ll  \lambda x + \mu.
	\end{align*}
	Equivalently, it suffices to prove for all $\lambda,\mu>0$
	\begin{align*}
	-\frac{\mu+1}{\lambda}\le-\frac{\mu}{\lambda},
	\end{align*}
	 which is obviously true.
Hence it follows that the matrix $G$ preserves interlacing polynomials.
Since $\left(f^{\langle r, r-1 \rangle}(x), \dots, f^{\langle  r,1 \rangle}(x), f^{\langle r,0 \rangle}(x)\right)$ is interlacing, so is $\big( g^{\langle r, r-1 \rangle}(x), \dots$, $g^{\langle  r,1 \rangle}(x), g^{\langle r,0 \rangle}(x) \big)$.
 This completes the proof.
\end{proof}

By iteratively using the above theorem, we obtain the following result. 

\begin{cor}\label{main2}
Let $r$ and $\ell$ be two  positive integers with $\ell \le r-1$.
 Suppose that
 \begin{align}\label{eq:h}
 \left(1+x+x^2+\dots+x^{\ell}\right)^n = h_{n,0}(x^r)+xh_{n,1}(x^r)+\dots+x^{r-1}h_{n,r-1}(x^{r}).
 \end{align}
 Then the polynomial sequence  $\left(h_{n,r-1}(x),\ldots, h_{n,1}(x), h_{n,0}(x)\right)$ is interlacing.
\end{cor}

\begin{proof}
We shall use induction on the integer $n$.
For the base $n=1$, the polynomial sequence $(0,\ldots, 0,1,\ldots,1)$ is interlacing by the convention that $0\ll  f$ for any real-rooted polynomial $f$.
Assume that the statement is true for $n=k$, that is to say that,  the polynomial sequence $\left(h_{k,r-1}(x),\ldots, h_{k,1}(x), h_{k,0}(x)\right)$ is interlacing. We proceed to show that $\left(h_{k+1,r-1}(x),\ldots, h_{k+1,1}(x), h_{k+1,0}(x)\right)$ is also interlacing. By \eqref{eq:h}, we have
\begin{align*}
 \left(1+x+x^2+\dots+x^{\ell}\right)\left( h_{k,0}(x^r)+xh_{k,1}(x^r)+\dots+x^{r-1}h_{k,r-1}(x^{r})\right) \\
 = h_{k+1,0}(x^r)+xh_{k+1,1}(x^r)+\dots+x^{r-1}h_{k+1,r-1}(x^{r}).
\end{align*}
Then, by the induction hypothesis and Theorem \ref{key}, we obtain the interlacing property of  the sequence $\left(h_{k+1,r-1}(x),\ldots, h_{k+1,1}(x), h_{k+1,0}(x)\right)$.
This completes the proof.
\end{proof}

Now we are in the position to prove Theorem \ref{main-thm}.

\begin{proof}[Proof of Theorem \ref{main-thm}]
Taking $\ell=r-2$  in \eqref{eq:h}, it follows that
 \begin{align*}
 \left(x+x^2+\dots+x^{r-2}+x^{r-1}\right)^n = x^nh_{n,0}(x^r)+x^{n+1}h_{n,1}(x^r)+\dots+x^{n+r-1}h_{n,r-1}(x^{r}).
 \end{align*}
Note that there exists one and only one integer in $\{0,1,\ldots, r-1\}$, say $j$, such that $r$ divides $n+j$.
Thus for this $j$,
\begin{align*}
\ {\rm E}_r \left(\, (x + x^2 + \cdots + x^{r-1})^n \right)=x^{(n+j)/r}h_{n,j}(x)
\end{align*}
by the definition of the linear operator ${\rm E}_r$. By Corollary \ref{main2}, the polynomial
$h_{n,j}(x)$ has only real zeros, so does $\ell_V \left((2^V)^{\langle r \rangle}, x \right)$.
This completes the proof.
\end{proof}

\vskip 3mm
\noindent {\bf Acknowledgments.}
 We would like to thank the referee for the valuable comments.
This work was supported by the National Science Foundation of China (Nos. 11626172, 11701424), the TJNU Funding for Scholars Studying Abroad, and  the PHD Program of TJNU (No. XB1616).


\end{document}